\documentclass[review, 11pt, authoryear]{elsarticle}

\usepackage{amsmath, amssymb, mathrsfs, amsthm, xcolor,orcidlink}
\usepackage[authoryear]{natbib}
\journal{Statistics \& Probability Letters}

\begin{document}

\begin{frontmatter}
\title{Bayesian Multiplicity Correction in the Probabilistic Forward Stepwise Framework}
\author[1]{Andrew Womack\corref{cor1}}
\ead{womack.andrew@gmail.com}
\author[2]{Daniel Taylor-Rodr\'{i}guez\fnref{fn2}}
\ead{dantayrod@pdx.edu}
\address[1]{Unaffiliated}
\address[2]{Department of Mathematics and Statistics, Portland State University}
\cortext[cor1]{Corresponding author}
\fntext[fn2]{Taylor-Rodr\'{i}guez was partially supported by NSF RTG DMS award 2136228.}
\date{April 1, 2026}

\begin{abstract}
We develop a natural Bayesian multiplicity-correcting prior distribution within the probabilistic forward stepwise representation of model space priors for regression problems. The proposed prior, obtained from making an analogy to the Holm procedure, exhibits behavior closely aligned with that of the Matryoshka doll prior. We compare both priors to several other priors, including some recently put forward as objective choices for model space prior probabilities. Our comparisons indicate that adequate multiplicity correction requires a degree of sparsity that many recommended priors do not provide, and we argue that multiplicity correction itself offers a principled and transparent criterion for specifying model space priors in regression.
%
\end{abstract}


\end{frontmatter}

\section{Introduction}
We are concerned with Bayesian multiplicity correction for model comparison in regression models with
\begin{equation}
\mathbb{E}[g(y_i)|\alpha, \boldsymbol{\beta}] = \alpha + \boldsymbol{\beta}^\top \boldsymbol{x}_i
\end{equation}
where $i=1,\ldots, n$ and $\boldsymbol{x}_i\in \mathbb{R}^p$. 
Let $[p]=\{1,\ldots,p\}$ and $\mathcal{P}_p=\{A\subseteq [p]\}$.
Bayesian multiplicity correction is the art of placing a prior distribution on $\mathcal{P}_p$ in such a way that the combinatorial number of subsets with $|A|=k$ is taken into account. This is equivalent to the specification of a prior distribution on model size $k=0, \ldots, p$. Let $\pi(k|p)=P(\{A:|A|=k\}$ and $\pi(k|\infty) = \lim_{p\uparrow\infty}\pi(k|p)$ be the limiting probabilities on model size as $p$ increases.

Specification of a prior on $\mathcal{P}_p$ is usually done using indicator variables $\gamma_j(A)$ for $j=1,\ldots, p$. Define $\gamma_j(A)$ to be $1$ if $j\in A$ and $0$ otherwise and define a distribution on $\{0,1\}^{\times p}$. The usual assumption is that the $\gamma_j$ are conditionally independent and identically distributed $\text{Bernoulli}(\phi)$ random variables and that $\phi$ follows some distribution. Making $\phi\perp\!\!\! \perp p$ leads to infinitely exchangeable priors on the $\gamma_j$ while letting the prior for $\phi$ depend on $p$ produces priors that are finitely exchangeable for a fixed $p$. A typical choice is a Beta prior on $\phi$ and thus a Beta-Binomial prior on the model space. 
\citet{ley2009effect} and \citet{scott2010bayes} recommend $\phi|p\sim\text{Beta}(1,1)$, which produces $\pi(k| p) = \frac{1}{p+1}$.
This produces a combinatorial accounting but does not provide any penalization on model size, leading to mass loss as $p$ increases. 
In contrast, \citet{wilson2010bayesian} recommend $\phi|p\sim\text{Beta}(1,\lambda p)$. The penalization on model size from this prior produces a Geometric$(\lambda/(1+\lambda))$ limiting distribution on model size as $p$ increases. 

Alternatively, one can specify a model space prior distribution by focusing on the poset structure of $\mathcal{P}_p$ directly. 
In \citet{wtf2025}, we developed the matryoshka doll (MD) prior following this approach.
The distribution is driven by the requirement that $P(A|p)$ be comparable to $P(\{B\in\mathcal{P}_p:B\supsetneq A\}|p)$ for all $A\in\mathcal{P}_p$.
This leads to a general condition that $\pi(k|p) = \omega_{p-k}\times(k+1)\times\pi(k+1|p)$ for some fixed sequence $(\omega_\ell)$. If we assume that $\omega_{\ell}\rightarrow \omega>0$ as $\ell$ increases, then the distribution on model size converges to a $\text{Poisson}(1/\omega)$ as $p$ increases. When we use the MD prior in the remainder of this paper, we will assume that the $\omega_\ell$ are constant (given by $\omega$), which directly produces a truncated Poisson on model size.

In this letter, we consider a third way
of building priors on $\mathcal{P}_p$, the probabilistic forward stepwise (pFS) representation of \citet{ma2015scalable}. The pFS specifies priors using paths of variable inclusion and so utilizes both indicator variable and poset structure considerations. 
We develop two natural multiplicity corrections in the pFS framework. 
These come from an analogy to the Holm procedure, which is a forward stepwise procedure. In Section \ref{sec:background} we review the pFS representation and the Holm procedure. In Section \ref{sec:priors} we develop the two new model space prior distributions and discuss their complexity penalization properties. 

Roughly, the two priors we obtain provide analogous behavior to the $\text{Beta-Binomial}(1,\lambda p)$ and the $\text{MD}(\omega)$ priors.
\cite{berger2026objectivemodelpriorprobabilities} discuss objective prior probabilities based on both indicator variables and direct assignment of $\pi(k| p)$. One prior they consider, the CMG prior \citep{casella_cluster_2014}
also strongly penalizes model size.
\cite{berger2026objectivemodelpriorprobabilities} argue against the objectivity of the $\text{Beta-Binomial}(1,\lambda p)$, $\text{MD}(\omega)$, and CMG priors calling them ``too parsimonious to be called objective.'' 
However, we contend that parsimony enforcement 
is a natural consideration when specifying an objective prior and these new pFS priors reinforce this notion. 


\section{Background}
\label{sec:background}

\subsection{The pFS Construction}
The pFS construction \citep{ma2015scalable} is based on probabilities for paths. A $k$-path in $[p]$ is a $k$-tuple (an element of $[p]^{\times k}$) without repeated elements. Let `$()$' represent the length $0$ empty path.
For each $k$-path, the pFS requires the specification of a $p-k+1$ dimensional conditional probability vector.
Let $\boldsymbol{i} = (i_1, \ldots, i_k)$ be a $k$-path and  
define the conditional probability vector $Q(\cdot|\boldsymbol{i},p)$ by defining conditional path probabilities $Q(j|\boldsymbol{i}, p)$ for each $j\in[p]\setminus\{i_1,\ldots,i_k\}$ as well as the conditional
path stopping probability $Q(\emptyset|\boldsymbol{i},p)$. 

The probability of the $k$-path $\boldsymbol{i}$ is given by
\begin{equation}
P(\boldsymbol{i}|p)
=
Q(i_1|(),p)
\cdot
Q(i_2|(i_1),p)
\times
\cdots\times
Q(i_k|(i_1,\ldots,i_{k-1}),p)
\cdot
Q(\emptyset|\boldsymbol{i},p)
\end{equation}
where the probability of the length $0$ path is $P(()|p) = Q(\emptyset|(),p)$.
A model is a subset $\{i_1,\ldots, i_k\}\subset[p]$ and can be reached through $k!$ different paths. 
Let $\mathcal{S}^k$ be the group of permutations on $\{1,\ldots,k\}$. The probability of the model $\{i_1,\ldots,i_k\}$ is 
\begin{equation}
P(\{i_1,\ldots, i_k\}|p)
=
\sum_{\sigma\in \mathcal{S}^k}
P((i_{\sigma(1)},\ldots, i_{\sigma(k)})|p).
\end{equation}

In order to get finite exchangeability, two restrictions must be made in the pFS construction. These restrictions are
\begin{align}
\label{eq:restriction1}
Q(j|\boldsymbol{i},p) &= \frac{1-Q(\emptyset|\boldsymbol{i},p)}{p-k}
\text{ for } j\in[p]\setminus\{i_1,\ldots,i_k\} \text{ and }
\\\label{eq:restriction2}
Q(\emptyset|\boldsymbol{i},p)
&=
Q(\emptyset|\boldsymbol{j},p)
\text{ for all }k\text{-paths }\boldsymbol{i}\text{ and }\boldsymbol{j}.
\end{align}
We define $Q_k(\emptyset|p)$ to be the common stopping probability for $k$-paths.
With these restrictions, the probability of a model (with $k>0$) is given by
\begin{equation}
P(\{i_1,\ldots, i_k\}|p)
=
\binom{p}{k}^{-1}
\times
\prod_{\ell=0}^{k-1}
\left(1-Q_\ell(\emptyset|p)\right)
\times
Q_k(\emptyset|p).
\end{equation}
This leads a prior on model size given by
\begin{equation}
\label{eq:symmetric}
\pi(k|p)
=
\prod_{\ell=0}^{k-1}
\left(1-Q_\ell(\emptyset|p)\right)
\times
Q_k(\emptyset|p)
\end{equation}
for $k=1,\ldots, p$ and $\pi_p(0)  = Q_0(\emptyset|p)$.

\subsection{The Holm Procedure}
The Holm procedure is a forward stepwise procedure designed to control the family-wise error rate, say at level $\alpha\in(0,1)$. Suppose there are $m$ hypotheses, $H_1,\ldots,H_m$, with sorted p-values $p_{(1)}\leq p_{(2)}\leq \cdots \leq p_{(m)}$. Let $H_{(1)}, \ldots, H_{(m)}$ be the hypotheses sorted following the ordering of their p-values. The Holm procedure rejects $H_{(1)}$ if $p_{(1)}<\frac{\alpha}{m}$. Otherwise, it fails to reject $H_{(1)}$ and the procedures stops. If the procedure is not stopped at $k-1$ steps, then $H_{(k)}$ is rejected if $p_{(k)}<\frac{\alpha}{m-k+1}$. Otherwise, it fails to reject $H_{(k)}$ and the procedures stops. If the Holm procedure stops at $k$ steps, the procedure rejects the $k-1$ hypotheses with smallest p-values and fails to reject the remaining $m-k+1$ hypotheses.

\section{pFS Holm Priors}
\label{sec:priors}

\subsection{Prior Definitions}

In analogy to the Holm procedure, when conditioned on a $k$-path  set \begin{equation}
Q(j|(i_1,\ldots,i_k),p) = \frac{\alpha}{p-k}
\end{equation}
for $j\in[p]\setminus\{i_1,\ldots,i_k\}$ and some $\alpha\in(0,1)$. This is equivalent to setting $Q_k(\emptyset|p) = 1-\alpha$ and provides
\begin{equation}
\label{eq:pathholm}
\pi(k|p)
=
(1-\alpha)\alpha^k
\end{equation}
for $k=0,\ldots, p-1$ and $\pi(p|p) = \alpha^p$. This prior is a truncated geometric distribution on model complexity (though the renormalization is achieved by increasing the probability assigned to $k=p$). This prior has the same behavior for large $p$ as a $\text{Beta-Binomial}(1,\lambda p)$ with $\lambda=\frac{1}{\alpha}-1$. 

The only issue with this multiplicity correction is that the Holm procedure is applied to paths and the model space prior forgets path ordering. One can correct the pFS Holm prior by designing a way to effectively divide $\pi(k|p)$ by the $k!$ multiplicity of $k$-paths. One solution is to divide out by the $k+1$ choices of insertion position of the next element that will be added to whatever $k$-path came before its inclusion.
This leads us to set
\begin{equation}
\label{eq:shp_def}
Q(j|(i_1,\ldots,i_k),p) = \frac{1}{k+1}\frac{\alpha}{p-k}
\end{equation}
for $j\in[p]\setminus\{i_1,\ldots, i_k\}$.
This is equivalent to setting $Q_k(\emptyset|p) = \frac{k+1-\alpha}{k+1}$ and we get
\begin{equation}
\label{eq:subsetholm}
\pi(k|p)
=
\frac{k+1-\alpha}{k+1}\times \frac{\alpha^k}{k!}
\end{equation}
for $k=0,\ldots, p-1$ and 
$\pi(p|p) = \frac{\alpha^p}{p!}$. 
This prior is quite similar to a truncated Poisson distribution and thus produces similar behavior to the matryoshka doll prior. 

A simple generalization that allows parameters that are greater than or equal to one is to set $Q_k(\emptyset|p) = \frac{k+\phi}{k+\phi+\theta}$ for $\phi, \theta > 0$. This provides
\begin{equation}
\label{eq:subsetholm_general}
\pi(k|p) 
= \frac{k+\phi}{k+\phi+\theta} \times \frac{\Gamma(\phi+\theta)\times \theta^k}{\Gamma(k+\phi+\theta)} 
\end{equation}
where $\Gamma(\cdot)$ is the gamma function and 
$\pi(p|p) = \frac{\Gamma(\phi+\theta) \times \theta^p}{\Gamma(p+\phi+\theta)}$. 
For $\alpha\in(0,1)$, the prior from \eqref{eq:subsetholm} is obtained by taking $\phi = 1-\alpha$ and $\theta=\alpha$. Alternatively, the prior from \eqref{eq:pathholm} is obtained by taking $\phi = \theta\times\frac{1-\alpha}{\alpha}$ and letting $\theta\uparrow\infty$.
For ease of exposition, we refer to the priors in \eqref{eq:pathholm}
and \eqref{eq:subsetholm_general} as the Path Holm Procedure (PHP) and Subset Holm Procedure (SHP) priors%
, respectively. The left panel in Figure \ref{fig:prior_considerations} shows the (logged) prior probability on model size for $p=20$ for the PHP and SHP priors as well as select other priors.

\subsection{Prior Effect on Learning}
%
The most important property for a model space prior, especially in the large $p$ regime, is how it penalizes the children set of a model. For $A\subsetneq[p]$, define the children set of $A$ by
\begin{equation}
\mathcal{C}_p(A) = \{B\in\mathcal{P}_p: B=A\cup\{j\}\text{ where }j\in [p]\setminus A\}.
\end{equation} 
Generally, the children set comprises the hardest models for a true model $A_T$ to beat in terms of Bayes Factors (with a learning rate that is, under mild conditions, $\sqrt{n}$ where $n$ is the sample size).

\begin{figure}[ht]
\centering
\includegraphics[width = 0.49\textwidth, page=6]{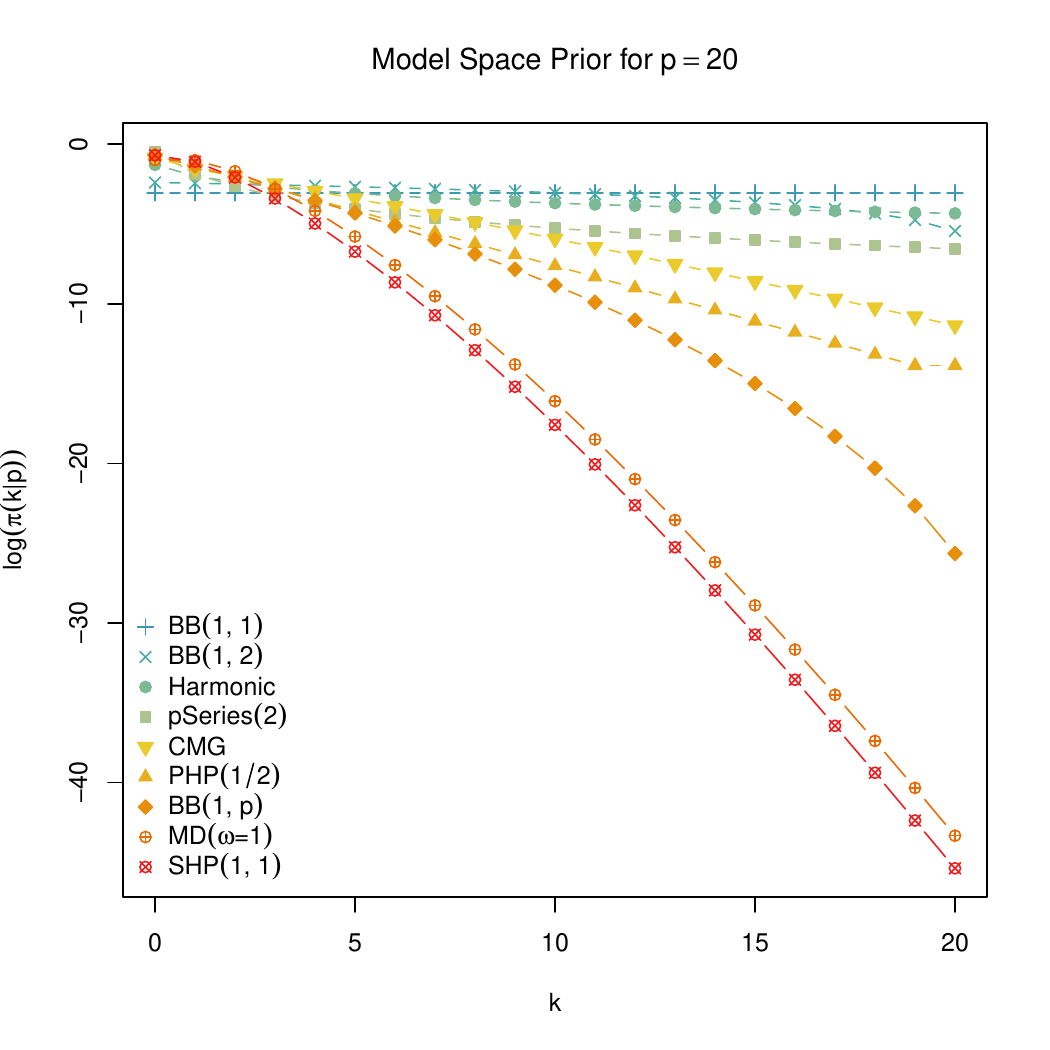}
\includegraphics[width = 0.49\textwidth, page=8]{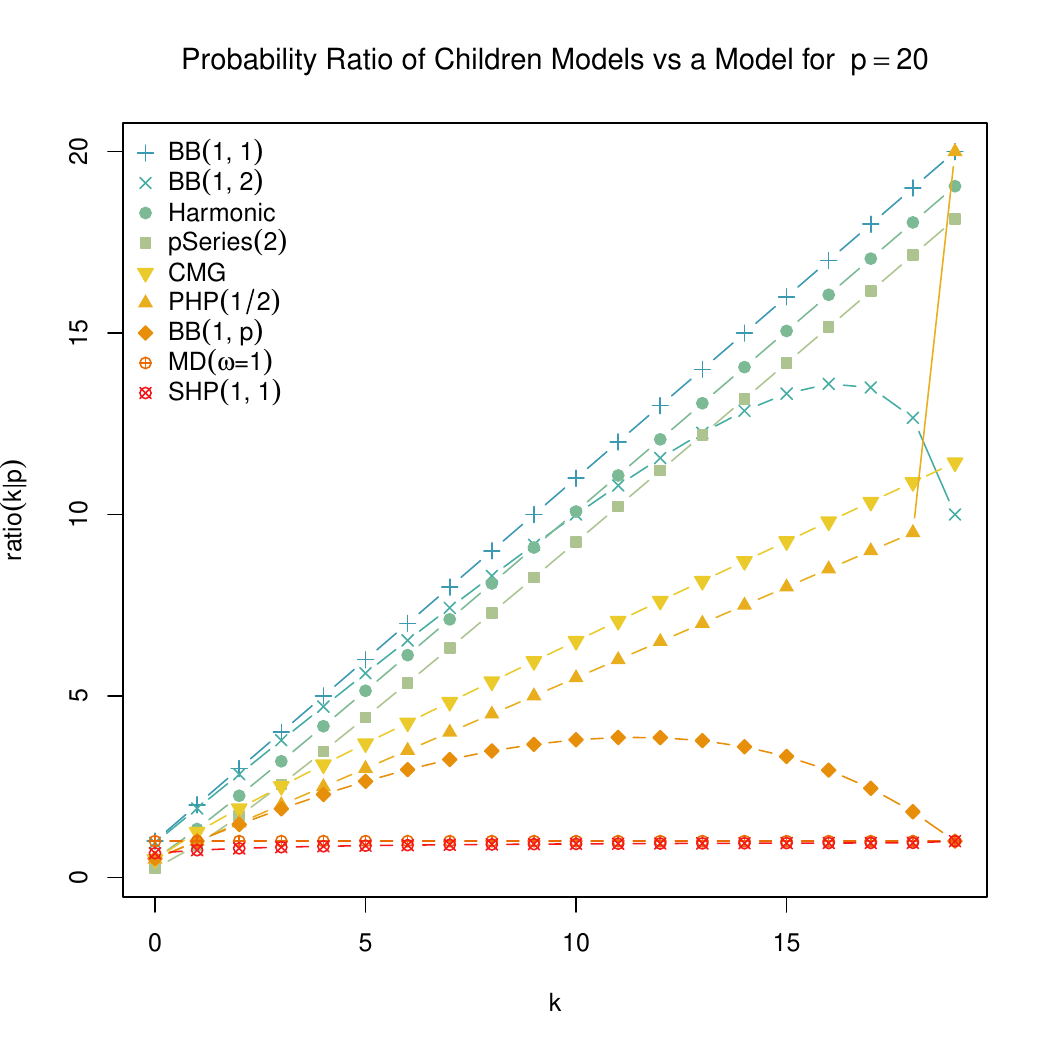}
\caption{Prior consideration of select priors. Left Panel: log-priors on model size. Right Panel: ratio of prior probavility of the set of children models versus a model (see \eqref{eq:rat}). For the Harmonic prior, $\pi(k|p) \propto (1+k)^{-1}$. For the P-Series prior, $\pi(k|p) \propto (1+k)^{-2}$.
For the CMG prior, $\pi(k|p)\propto E[Y^{2k}]/k!$ where $Y\sim\mathcal{N}(0.5, 0.25)$.}
\label{fig:prior_considerations}
\end{figure}

The ratio of prior probabilities of the children set of $A$ versus $A$ (assuming $k=|A|<p$) is given by
\begin{equation}
\label{eq:rat}
\text{ratio}(k|p) = 
\frac{P(\mathcal{C}_p(A)|p)}{P(A|p)}
=
\frac{(k+1)\times\pi(k+1| p)}{\pi(k |p)}.
\end{equation}
For the PHP and SHP priors (assuming $0\leq k<p-1$) we get
\begin{align}
\text{ratio}_{\text{PHP}}(k|p) 
&= 
(k+1)\times\alpha
\\
\text{ratio}_{\text{SHP}}(k|p) 
&
= 
\frac{(k+1+\phi)(k+1)}{(k+1+\phi+\theta)(k+\phi)}
\times
\theta.
\end{align}
The difference is stark, the PHP prior has $\text{ratio}(k|p)$ increasing linearly with $k$ whereas the SHP prior has $\text{ratio}(k|p)$ converging to $\theta$ as $k$ increases. The behavior of the PHP prior is not really problematic when the size of the true model, $p_T = |A_T|$, is fixed and the sample size $n$ is large. However, once $p_T\propto \sqrt{n}$, then the prior behavior of the PHP prior effectively cancels out the Bayes Factor learning rate and it becomes impossible for the posterior probability of $A_T$ to converge to one. Even if $p_T = o(\sqrt{n})$ as $p_T$ increases with $n$, the learning rate against the set of children models decreases as $n$ increases for the PHP prior. The SHP prior
does not 
exhibit this bad behavior regardless of the growth rate of $p_T$ with $n$. 
The right panel in Figure \ref{fig:prior_considerations} shows \eqref{eq:rat} for these priors as well as select other priors which will be included in the simulation study in Section \ref{sec:simulation}.

\clearpage

\section{Simulation}
\label{sec:simulation}
We perform two simulations with large $p$ in order to demonstrate the differences in the priors. In addition to the PHP and SHP priors, we also include the other priors from Figure \ref{fig:prior_considerations}. These choices cover the priors of \cite{scott2010bayes}, \cite{wilson2010bayesian}, \cite{casella_cluster_2014}, and \cite{wtf2025} as well as two others discussed in \cite{berger2026objectivemodelpriorprobabilities}. The The P-Series$(\text{power}=2)$ is included to provide a prior on model complexity that does not exhibit mass loss as $p\rightarrow\infty$ and has the same asymptotic penalization profile as the Harmonic prior and fixed Beta-Binomial priors.

The first set of simulations uses $n=200$, $p\in\{n/2, n, 3n/2\}$, $p_{T}\in\{5, 10, 15\}$, signal to noise ratio $4$. The second set of simulations takes $n\in\{125, 250, 500\}$, $p=2n$, $p_T = p/25$, and signal to noise ratio in $\{4, 5.67, 9\}$ (corresponding to $R^2\in\{0.8, 0.85, 0.9\}$). Both simulations use  independent $\text{Normal}(0,1)$ covariates, the Zellner-Siow prior \citep{zellner1980posterior} for Bayes Factor computation, and the Moore-Penrose pseudo-inverse to define a singular normal prior on regression coefficients whenever a model's design matrix is not full column rank. Models were fit using MCMC ($10^6$ draws) with single indicator complementation, in-out swap, and same-size model replacement proposal kernels.

\begin{figure}[ht]
\centering
\includegraphics[width = 0.49\textwidth]{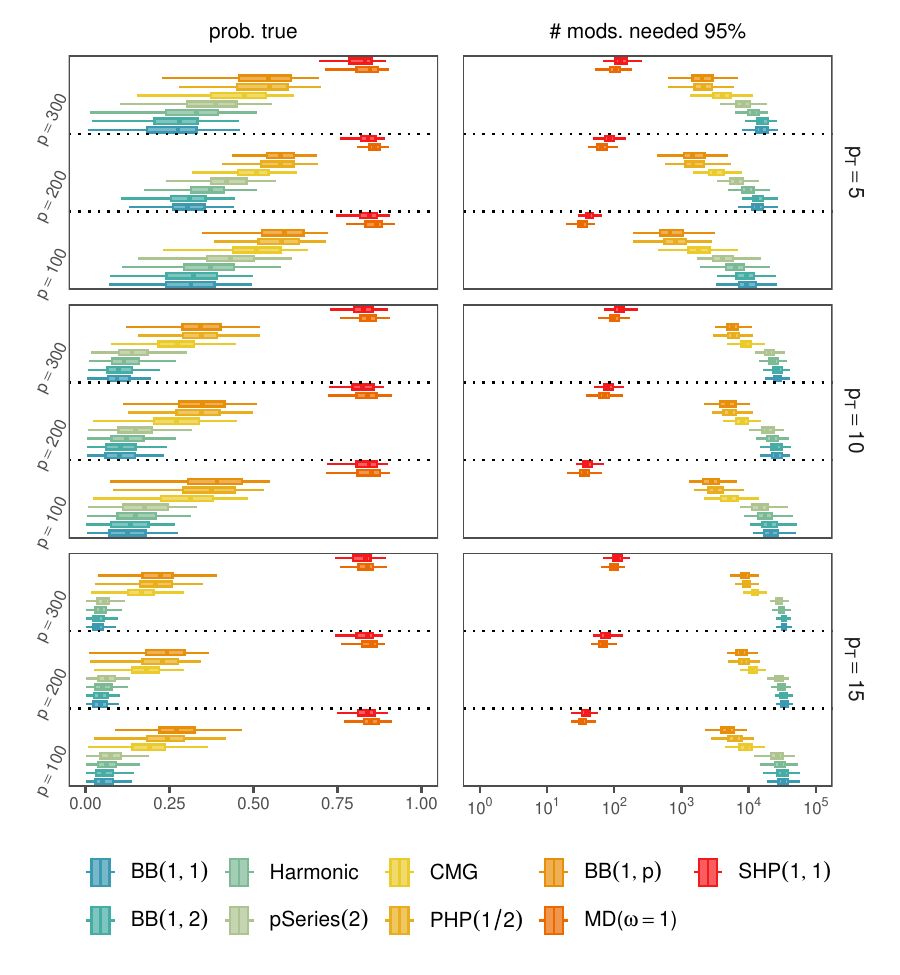}
\includegraphics[width = 0.49\textwidth]{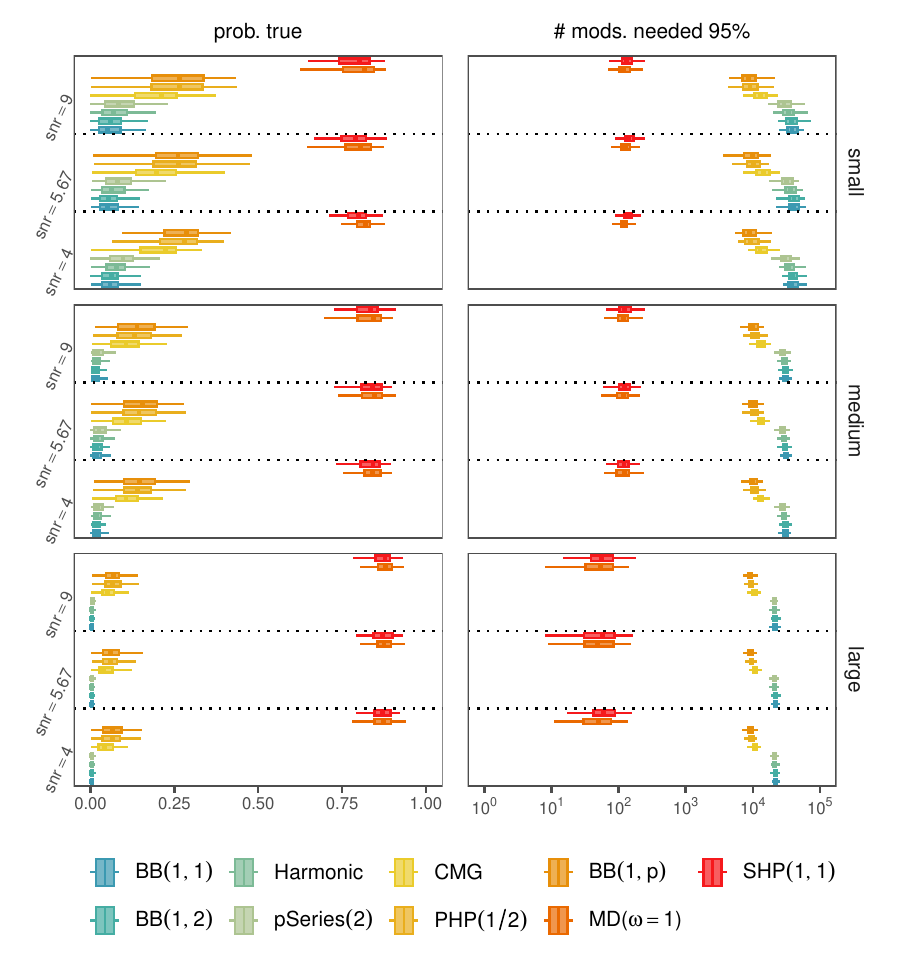}
\caption{Simulation results from 100 replications with $10^6$ MCMC draws per replication. Plotted measurements are the posterior probability of the true model and the number of models needed to obtain $95\%$ of the posterior mass. All probabilities are assessed using sampling frequency. Left Panel: $p$ and $p_T$ given on the left and right axis labels with fixed $n=200$ and signal to noise ratio $4$. Right Panel: signal to noise ratio given on the left axis label and simulation scenario given on the right axis label; for this scenario $p=2n$ and $p_T=p/25$ with $n\in\{125, 250, 500\}$ labeled as small, medium, and large. }
\label{fig:simulation}
\end{figure}

Figure \ref{fig:simulation} shows the results of 100 replicates of the simulations. The pattern is clear, Poisson-like penalization on model size is necessary for posterior probability to concentrate on the true model in the large-$p$ regime. The separation of the models into three distinct groups is clearly seen from the number of models necessary to achieve $95\%$ of the posterior probability. The $\text{Beta-Binomial}(1,1)$, $\text{Beta-Binomial}(1,2)$, Harmonic, and  $P-Series(\text{power}=2)$ priors for the first group. The priors provide $\text{ratio}(k|p)\approx k+1$ as $p\rightarrow\infty$ and posterior probability simply cannot concentrate on the true model. The second group is comprised of the $\text{PHP}(\alpha=0.5)$, $\text{Beta-Binomial}(1,p)$, and CMG  priors. The value for $\text{ratio}(k|p)$ is $(k+1)/2$ as $p\rightarrow\infty$ for the $\text{PHP}(\alpha=0.5)$ and $\text{Beta-Binomial}(1,p)$ priors. 
For the CMG, it is easy to show that the ratio is the open interval $((k+1)/2, k+1)$ and close to $(k+1)/2$, but an exact form eludes us. Though providing better posterior inference than the first set or priors, this second group is similarly incapable of concentrating posterior probability on the true model. The third group is the SHP$(\phi=\theta=1)$ and MD$(\omega=1)$ priors. These priors provide Poisson-like penalization on model size and the ratio is $1$ (exactly for the MD and approximately for the SHP). These are the only two priors under consideration that are capable of concentrating posterior probability on the true model, especially in the regime where $p$ and $p_T$ grow with $n$.

\section{Discussion}
\label{sec:discussion}


Multiplicity correction, especially in high dimensional problems, provides particular difficulties in the Bayesian paradigm. 
We posit that this is because of the standard way of assigning priors using exchangeable or conditionally independent inclusion indicators. 
Neither the fixed $\text{Beta-Binomial}(a,b)$ nor the $\text{Beta-Binomial}(a,b+\lambda p)$ priors (or priors that behave like them in terms of \eqref{eq:rat}) sufficiently penalize model size to control false positive inclusion in the posterior.
Multiplicity correction needs to be explicitly built into the prior construction in order to obtain appropriate control on false discovery inclusion and produce posterior concentration on the true model.

The Subset Holm Procedure prior multiplicity-corrects the path structure of the model space,
while the Matryoshka Doll prior does so through poset considerations. The fact that these two rather different constructions converge on the same Poisson-like behavior is, we think, one of the most interesting findings of this work. It suggests something that is intrinsic to the geometry of standard regression model spaces rather than an artifact of either construction.

These priors are objective in the sense that they are
derived by specifying a general rule, procedure, or desirable property that is then operationalized in this specific context. One can apply the same rule to more structured model spaces and obtain default
prior distributions in those settings. 
The advantage of using multiplicity correction as the 
organizing principle is that it acknowledges the structure of the model space as it pertains to inference. The weak Bayes Factors learning rate against false positives (versus false negatives) requires the model space prior to place disproportionate mass on parsimonious models. Diffuse priors do not sufficiently promote parsimony and produce posterior distributions that cannot concentrate on the true model. Proper accounting for multiplicity yields priors that strongly enforce parsimony and produce appropriate posterior behavior, especially in the large-$p$ regime.

\bibliographystyle{elsarticle-harv} 
\bibliography{wtf-2}

\end{document}